\documentclass{article}

\usepackage{graphicx}
\usepackage{hyperref}
\usepackage{amsmath}
\usepackage{amssymb}

\numberwithin{equation}{section}

\newcommand{\D}{{\mathrm D}}
\newcommand{\F}{{\mathcal{F}}}

\newcommand{\dd}{\mathrm{d}}

\newcommand{\e}{\mathrm{e}}

\newcommand{\R}{\mathbb{R}}

\newcommand{\Sc}{{\mathcal{S}}}
\newcommand{\Hs}{\dot H}

\newtheorem{thm}{Theorem}
\newtheorem{prop}{Proposition}
\newtheorem{coro}[thm]{Corollary}
\newtheorem{defi}{Definition.}
\newtheorem{lem}[defi]{Lemma.}
\newtheorem{rem}{Remark}

\newenvironment{prf}{\noindent {\bf Proof} }{\endprf\par}
\def \endprf{\hfill  {\vrule height6pt width6pt depth0pt}\medskip}

\title{Some Remarks on Strichartz Estimates for Homogeneous Wave Equation \thanks {This work is supported by
NSFC 10271108. }}
\author{Daoyuan Fang and Chengbo Wang\thanks {email: DF:
dyf@zju.edu.cn, CW: wangcbo@yahoo.com.cn}
\\
Department of Mathematics, Zhejiang University,\\
Hangzhou 310027, P. R. China} \date{}
\begin{document}
\maketitle

\begin{abstract}
We give several remarks on Strichartz estimates for homogeneous
wave equation with special attention to the cases of $L^\infty_x$
estimates, radial solutions and initial data from the
inhomogeneous Sobolev spaces. In particular, we give the failure
of the endpoint estimate $L^4_t L^\infty_x$ for $n=2$.
\end{abstract}

\section{Introduction}
As is well-known, Strichartz-type estimate is of particular
importance in the low regularity well-posedness theory for
semilinear wave and Schr\"{o}dinger equations, e.g. in
\cite{PoSi93}, \cite{LdSo95} and \cite{CaWe90}. Recently, there
are many diverse advances in extending Strichartz-type estimates
for both wave equations and Schr\"{o}dinger equations, such as in
\cite{Tao00}, \cite{Fc03p}, \cite{MaNaNaOz03p}, \cite{Stbz04p}
etc. Since the appearance of \cite{KeTa98}, it is generally
believed that the Strichartz estimates of homogeneous equation has
been totally solved and the only remained problem is to extend it
to the inhomogeneous equation. However, as we know, there are
still some gaps for homogeneous estimate of wave equation and it
seems that it causes some confusions(many authors have different
statement concerning such estimate).

In this paper, we'll concern solely on Strichartz estimate and its
variants for homogeneous wave equation. For the counterpart of
Schr\"{o}dinger equation, one may consult \cite{KeTa98} and
\cite{Tao00}. As usual, we denote the space of Schwartz class by
$\Sc$ and use $\hat{f}=\F(f)$ denote the Fourier transform of
$f\in\Sc'$ and let $p(\D)f=\F^{-1}(p(|\xi|) \hat{f}(\xi))$. Also,
we use $\Hs^s$ to denote the usual homogeneous Sobolev space
$\D^{-s}L^2(\mathbb{R}^n)$ for $s<n/2$(note that for $s\geq n/2$,
one would interpret such spaces as the subspace of $\Sc'$ modulo
polynomials of degree less than or equals $[s-n/2]$). Moreover, we
define the homogeneous Besov space $\dot{B}^{s}_{p,q}$ for $s<n/p$
or $s=n/p$ with $q=1$ as follows. Let $\triangle_j f
:=\F^{-1}(\varphi(2^{-j}\xi)\hat{f}(\xi))$ be the usual
homogeneous Littlewood-Paley projection,
$\|f\|_{\dot{B}^{s}_{p,q}}=\|2^{j s} \triangle_j f\|_{l^q_j
L_x^p}$. Then $\dot{B}^{s}_{p,q}=\{f\in\Sc',\
\|f\|_{\dot{B}^{s}_{p,q}}<\infty,\ {\rm and}\ \sum_{j\in
\mathbb{Z}}\triangle_j f=f\ {\rm in}\ \Sc'\}$. For general $s$,
one should introduce such space in $\Sc'$ modulo finite degree
polynomials. Note also that $u = \cos(t \D) f+\D^{-1}\sin(t \D) g$
solves homogeneous wave equation $\square u=0$ with data $(f,g)$,
so we only need to state the estimate for operator $\exp(i t \D)$.

At first, we give a definition.
\begin{defi} Let $n\geq 2$ and $2\leq q,r\leq \infty$, we say that the triple $(q,r,n)$ is admissible if
\begin{equation} \label{admi}\frac{1}{q}\leq
\frac{n-1}{2}(\frac{1}{2}-\frac{1}{r}).
\end{equation} And we say the triple is radial-admissible if \eqref{admi} is
substituted by \begin{equation} \label{radial}\frac{1}{q}<
(n-1)(\frac{1}{2}-\frac{1}{r})\ .\end{equation}
\end{defi}

The classical Strichartz-type estimates are essentially the
following single frequency estimate:
\begin{thm}{\rm (Essential Strichartz Estimate)}\label{ess}
Let $n\geq 2$, then the following two statements are equivalent,\\
(I) the single frequency estimate \begin{equation}
\label{Stri}\|\exp(i t \D) f(x)\|_{L^q_t L_x^r}\lesssim \|
f\|_{L_x^2}\ ,\end{equation} valid for all
$f\in L^2$ with $supp(\hat{f})\subset\{1/2<|\xi|<2\}$;\\
(II) (q,r,n) is admissible and $(q,r,n)\neq(2,\infty,3)$.
\end{thm}
The positive results are given in \cite{GiVe95} and \cite{KeTa98}.
The necessary condition $q\geq 2$ is given by time-translation
invariant argument, and \eqref{admi} follows from A.W.Knapp's
counterexample $\hat{f}=\chi_A$ with $A=\{\xi | 1/2<\xi_1<3/2,
|\xi_j|<\epsilon, 2\leq j\leq n \}$ by letting
$\epsilon\rightarrow 0$(note that this example is non-radial). The
forbidden triple $(q,r,n)=(2,\infty,3)$ is given in
\cite{Tao}(with previous results in \cite{KlMa93} and
\cite{Mo98}).

By applying homogeneous Littlewood-Paley decomposition and
scaling, one would get immediately the following(throughout this
paper, $b=n(\frac{1}{2}-\frac{1}{r})-\frac{1}{q}$)
\begin{coro}{\rm (Classical
Strichartz Estimate)}\label{cls}Let
$b=n(\frac{1}{2}-\frac{1}{r})-\frac{1}{q}$, then for all
admissible $(q,r,n)$ with $r<\infty$, \begin{equation}
\label{Sobolev}\|\exp(i t \D) f(x)\|_{L^q_t L_x^r}\lesssim \|
f\|_{\dot{H}^b}\ ,\end{equation} Moreover, for all admissible
$(q,r,n)$ except that $(q,r,n)=(2,\infty,3)$, \begin{equation}
\label{Besov}\|\exp(i t \D)
f(x)\|_{L^q_t\dot{\mathbf{B}}^0_{r,2}}\lesssim \| f\|_{\dot{H}^b}\
,\end{equation} \begin{equation} \label{Besov2}\|\exp(i t \D)
f(x)\|_{L^q_tL^r_x}\lesssim \| f\|_{\dot{\mathbf{B}}^b_{2,1}}\
.\end{equation}\end{coro}

Then a question arises naturally:{\sf Can \eqref{Sobolev} valid
with $r=\infty$? Or if it fails, how can it be improved to
restricted cases such as spherically symmetry or angular
regularity?}

In fact, we have the following results by supplement the known
results.
\begin{thm}\label{gen}
Let $b=n(\frac{1}{2}-\frac{1}{r})-\frac{1}{q}$, then we have
\eqref{Sobolev} for all admissible $(q,r,n)$ except that
$(q,r)=(\max(2,\frac{4}{n-1}),\infty)$ and
$(q,r)=(\infty,\infty)$. On the other hand, in order for
\eqref{Sobolev} to be valid for all $f\in\Sc$, we need $(q,r,n)$
admissible, $(q,r)\neq(\infty,\infty)$ and $q>\frac{4}{n-1}$.
\end{thm}

\begin{figure}\centering
  \includegraphics[width=0.90\textwidth]{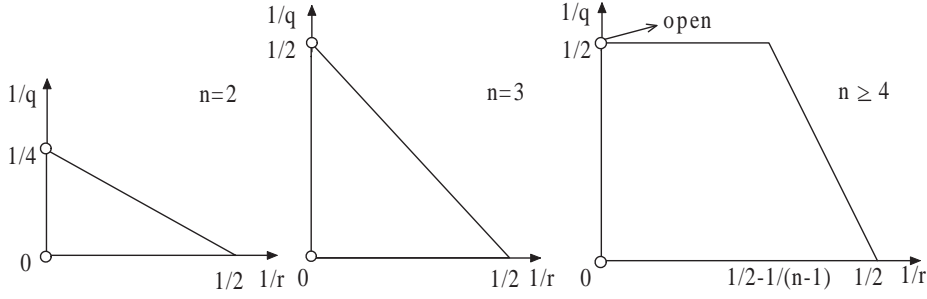}\\
  \caption{Classical Strichartz Estimate \eqref{Sobolev}}\label{cls-p}
\end{figure}

\begin{rem}
As stated in Theorem \ref{gen}, the only remained open problem for
homogeneous estimate now is the endpoint $(2,\infty,n)$ with
$n\geq 4$. Alternatively, we have(ref Proposition \ref{substi})
$$\|\exp(i t \D)
f(x)\|_{L^2_t L_x^{\infty}}\lesssim \|
f\|_{\dot{H}^{b-\epsilon}}^{\theta}\|
f\|_{\dot{H}^{b+\delta}}^{1-\theta}\ .$$ For the valid region of
$(1/q,1/r)$ for \eqref{Sobolev}, see Figure \ref{cls-p}.
\end{rem}
\begin{rem}
For  $r=\infty$, it seems that there are some confusions for the
validness of estimate \eqref{Sobolev}, except for $(2,\infty,3)$
and $(\infty,\infty,n)$. Many authors have different statements
for such estimate. So the partial purpose of Theorem \ref{gen} is
to clarify this confusion.
\end{rem}

\begin{rem}
Note that the embeddings $\dot{H}^{n/2}\subset L^\infty$ and
$H^{n/2}\subset L^\infty$ are both fail to valid even for radial
function, then \eqref{Sobolev} fails for $(q,r)=(\infty,\infty)$
for any $n$.
\end{rem}

\begin{rem}As a complement to the failure of some $r=\infty$ estimate in
\eqref{Sobolev}, we have a simple but somewhat interesting
result(Proposition \ref{anti}). Let $2\leq q<\infty$,
we have
$$\|\exp(i t \D)f(x)\|_{L_x^\infty L_t^q}
\lesssim \|f\|_{\dot H^b}\ .$$
\end{rem}

For radial function,  the region of "admissible" triple can be
vastly improved(the angular improvement will be given in Theorem
\ref{Ster} of Section 3).
\begin{thm}\label{radi}
Let $(q,r,n)$ be radial-admissible and $(q,r)\neq(\infty,\infty)$,
then \eqref{Sobolev} valid for all radial $f$.
\end{thm}
\begin{figure}
\centering
  \includegraphics[width=0.70\textwidth]{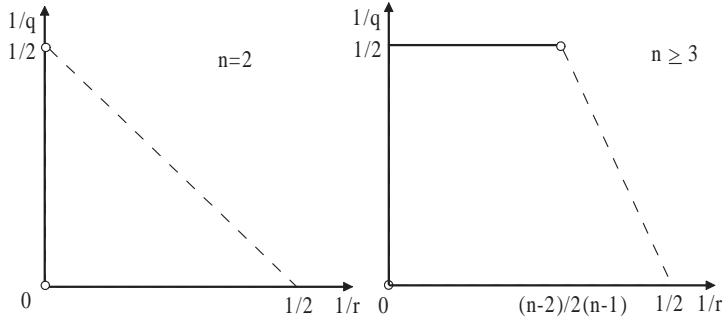}\\
  \caption{Radial improvement of Strichartz Estimate \eqref{Sobolev}}\label{radial-p}
\end{figure}

For the completeness of exposition, we also state here the
correspond estimate in the Sobolev space $H^s$. Here $b+$ denotes
the $b+\epsilon$ with $\epsilon>0$ arbitrary small.
\begin{thm}\label{inho} Let $n\geq 2$, and $u(t,x)$ be the solution to $\square u=0$ with data
$(f,g)$. Then the estimate \begin{equation}
\label{inh}\|u\|_{L_t^q L_x^r} \lesssim
\|f\|_{H^s}+\|g\|_{H^{s-1}}\end{equation} valid with $s=b+$ 
if $n\geq 3$, $b\geq 1$, $(q,r,n)$ admissible and $(q,r,n)\neq
(2,\infty,3)$,
on the other
hand, \eqref{inh} valid with $s=b+$ only if 
$b\geq 1$, $(q,r,n)$ admissible and $(q,r,n)\neq (2,\infty,3)$.
Moreover, if $n\geq 3$, $b\geq 1$, $(q,r,n)$ admissible and
$(q,r)\neq(2,\infty)$, $(\infty,\infty)$, then \eqref{inh} valid
with $s=b$. And if \eqref{inh} valid with $s=b$, we need $n\geq
3$, $b\geq 1$, $(q,r,n)$ admissible, $(q,r,n)\neq(2,\infty,3)$ and
$(q,r)\neq(\infty,\infty)$.
\end{thm}
Note that \eqref{inh} valid for $s$ implies it's validness for
$s+$, also the $s+$ failure implies the $s$ failure.
\begin{figure}\centering
  \includegraphics[width=1.0\textwidth]{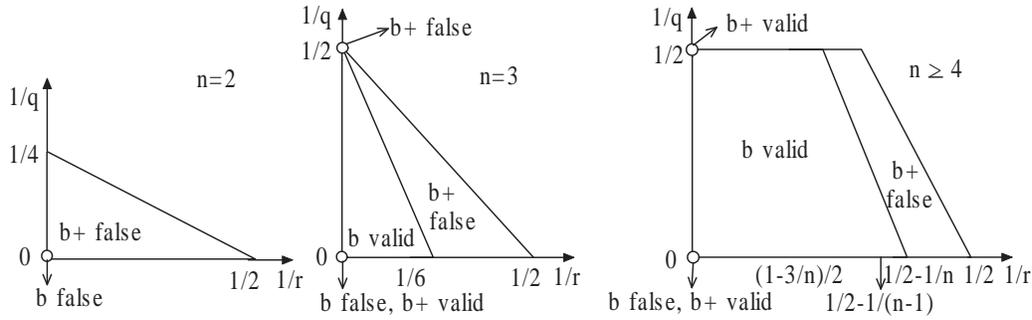}\\
  \caption{Theorem \ref{inho}: Strichartz Estimate in $H^{s}$
  }\label{inh-p}
\end{figure}

From the figure \ref{inh-p} of Theorem \ref{inho}, we see that
there is a new limitation for $H^s$ estimate, due to the fact that
$b-1$ may less than $0$. However, if one substitute $d u$ for $u$
in \eqref{inh} with $d u=(\partial_t u, \nabla u)$, one can
eliminate
such limitation. In fact, 
$$\|\exp(i t \D) f \|_{L_t^q L_x^r}
\lesssim \|f\|_{H^s}\ \Rightarrow\ \|d u\|_{L_t^q L_x^r} \lesssim
\|f\|_{H^{s+1}}+\|g\|_{H^{s}}\ \Rightarrow\ \|\cos(t \D) f
\|_{L_t^q
L_x^r}\lesssim \|f\|_{H^s}$$ 
Hence
\begin{thm}\label{d-inh}
Let $n\geq 2$, and $u(t,x)$ be the solution to $\square u=0$ with
data $(f,g)$. Then we have \begin{equation} \label{drv}\|d
u\|_{L_t^q L_x^r}\lesssim
\|f\|_{H^{s+1}}+\|g\|_{H^{s}}\end{equation} with $s=b+$ if and
only if $(q,r,n)$ admissible and $(q,r,n)\neq (2,\infty,3)$.
Moreover, for admissible $(q,r,n)$ except that
$(q,r)=(\max(2,\frac{4}{n-1}),\infty)$ and
$(q,r)=(\infty,\infty)$, we have \eqref{drv} valid with $s=b$. On
the other hand, if \eqref{drv} valid with $s=b$, then $(q,r,n)$ is
admissible and $(q,r,n)\neq (2,\infty,3),(\infty,\infty,n)$.
\end{thm}
\begin{figure}\centering
  \includegraphics[width=0.90\textwidth]{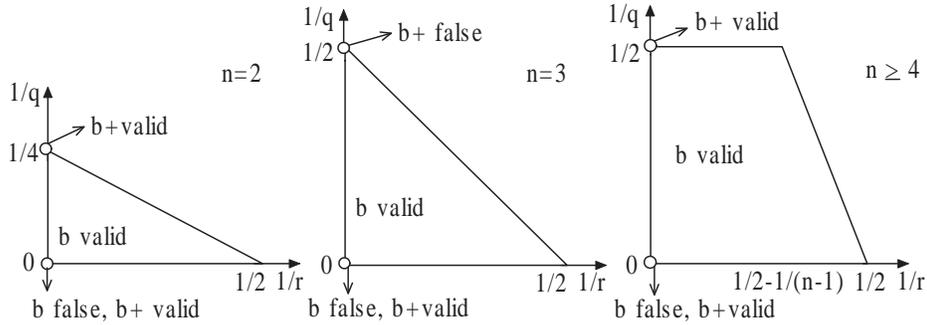}\\
  \caption{Theorem \ref{d-inh}, Strichartz Estimate in $H^{s}$
  }\label{drv-p}
\end{figure}

\section{General Case with $r=\infty$}

In this section, we concern on Theorem \ref{gen}. Firstly, we give
some remarks on a particular triple. The forbidden triple
$(q,r,n)=(2,\infty,3)$ achieves special attention in the study of
Strichartz estimate. \begin{itemize} \item This forbidden triple
for \eqref{Sobolev}is found in \cite{KlMa93}, they also find that
the triple is "admissible" if the data is radial symmetric.

\item In \cite{Mo98}, the author gives this triple's inadmissible
with the $L^\infty$ norm substitute by $BMO$ in \eqref{Sobolev}.
The correspond assertion in Theorem \ref{ess} may be found in
\cite{Tao}, which implies the result in \cite{Mo98}.(Since the
Littlewood-Paley projection of $BMO$ is in $L^\infty$).

\item However recently, in \cite{MaNaNaOz03p}, the authors show
that if one substitutes the $L^2L^\infty$ norm by $L^2_t
L^\infty_rL^p_\theta$ for any $p<\infty$, then \eqref{Sobolev}
will be valid with the operator bound like $\sqrt{p}$ as
$p\rightarrow\infty$. (Instead, as noted there, one may get the
$L^2L^\infty$ bound for the data with some angular regularity by
Sobolev embedding $H^{\epsilon,p}_\theta\subset L^\infty$, $p
\epsilon>2$. This fact will be used in Theorem \ref{Ster} below.)
\end{itemize}

As stated in Theorem \ref{gen}, there are many $q$ such that
\eqref{Sobolev} valid with $r=\infty$. In fact, Klainerman and
Machedon has stated this and explicitly prove the $n=3$ case in
\cite{KlMa98}. We restate it here.
\begin{prop}\label{gen-v}
If $\max(\frac{4}{n-1},2)<q<\infty$ and $r=\infty$, then
\eqref{Sobolev} valid.
\end{prop}
The proof of this result in \cite{KlMa98} is subtle. However it
follows easily from classical Strichartz estimate in Corollary
\ref{cls} and the following generalized Gagliardo-Nirenberg
inequality (Lemma 1.4 in \cite{EsVe97})
\begin{lem}\label{GN} Let $a,c\in(1,\infty)$, $\alpha,\beta\in(0,n)$ and $\alpha
a < n < \beta c$, then
$$\|f\|_{L^\infty}\lesssim \|\D^\alpha f\|_{L^a}^\theta \|\D^\beta
f\|_{L^c}^{1-\theta}\ ,$$ where
$\theta=(\frac{\beta}{n}-\frac{1}{c})/(\frac{\beta}{n}-\frac{1}{c}+\frac{1}{a}-\frac{\alpha}{n})$.
\end{lem}
Hence for $\max(\frac{4}{n-1},2)<q<\infty$, let $a=2$,
$\alpha=b<n/2$, one may choose $2\leq c<\infty$ and $\beta$ such
that $((1-\theta)q,c,n)$ admissible. In fact, one only needs to
choose
$\beta=\min(\frac{n}{c}+\frac{1}{2}-\frac{1}{q},\frac{n+1}{2
c}+\frac{n-1}{4}-\frac{1}{q})>\frac{n}{c}$. Then we can apply
Lemma \ref{GN} to yield \begin{eqnarray*} \|\exp(i t \D)f\|_{L_t^q
L_x^\infty}&\lesssim &(\|\D^b \exp(i t \D)f\|_{L_x^2}^\theta
\|\D^\beta \exp(i t \D)f\|_{L_x^c}^{1-\theta})_{L_t^q}\\
&\lesssim& \|f\|_{\Hs^b}^\theta \|\D^\beta \exp(i t
\D)f\|_{L_t^{q(1-\theta)}L_x^c}^{1-\theta}\lesssim \|f\|_{\Hs^b}\
.\end{eqnarray*} So is Proposition \ref{gen-v}. Similar argument
yields the following
\begin{prop}\label{substi}
Let $n\geq 4$ and $b=\frac{n-1}{2}$, then  $$\|\exp(i t \D)
f(x)\|_{L^2_t L_x^{\infty}}\lesssim \|
f\|_{\dot{H}^{b-\epsilon}}^{\theta}\|
f\|_{\dot{H}^{b+\delta}}^{1-\theta}$$ with $\epsilon\in
(0,\frac{n(n-3)}{2(n-1)})$, $\delta\in(0,n)$ and
$\delta=(\delta+\epsilon)\theta$.
\end{prop}

Now let's restate the stated necessary condition here.
\begin{prop}\label{gen-f}
If \eqref{Sobolev} valid, then $q>\frac{4}{n-1}$.
\end{prop}
\begin{prf}
Note that Knapp's counterexamples and time-translation invariant
argument give a necessary condition of $\frac{1}{q}\leq
\min(\frac{1}{2},\frac{n-1}{2}(\frac{1}{2}-\frac{1}{r}))$, then
the only remained triple is $(\frac{4}{n-1},\infty,n)$ with
$n=2,3$. We use the contradiction argument. Assume that
\eqref{Sobolev} valid for such triple, then for any $g\in \Sc(\R)$
and $f\in\Sc(\R^n)$, we have
$$|\int g(t) \langle \delta(x-t
\e_1), \D^{-b}\e^{it\D} f(x) \rangle_{x}\dd t| \leq
\|{\left(\D^{-b}\e^{it\D} f\right)(t,t\e_1)}\|_{L^q_t}
\|g\|_{L^{q'}_t} \lesssim \|f\|_{L^2}\|g\|_{L^{q'}_t}$$ with
$\e_1=(1,0,\cdots,0)$. Apply Plancherel theorem to the left hand
side with respect to $x$ yield that
$$\langle \hat{f}(-\xi),\int |{\xi}|^{-b} \e^{-i
t(|{\xi}|+\xi_1)}g(t)\dd t\rangle_\xi \lesssim
\|{f}\|_{L^2}\|{g}\|_{L^{q'}_t}\ ,$$

$$\langle \hat{f}(-\xi),|{\xi}|^{-b} \hat{g}(|{\xi}|+\xi_1)\rangle_\xi \lesssim
\|{f}\|_{L^2}\|{g}\|_{L^{q'}_t}\ .$$ From this, one has
\begin{equation} \label{I}\|{|{\xi}|^{-b}
\hat{g}(|{\xi}|+\xi_1)}\|_{L_{\xi}^2}\lesssim \|{g}\|_{L^{q'}_t}\
.\end{equation} Make coordinate transformation (denote
$\xi=(\xi_1,\xi')$ and $\lambda=(\lambda_1,\lambda')$) $
\lambda_1=|{\xi}|+\xi_1\geq 0$, $\lambda'=\xi'$, we have
$\xi_1=\frac{\lambda_1^2-|{\lambda'}|^2}{2\lambda_1}$,
$|{\xi}|=\frac{|{\lambda}|^2}{2\lambda_1}$, and $\dd
\xi=\frac{|{\lambda}|^2}{2\lambda_1^2}\dd \lambda$. Note that
$2(1-2b)=1-n$ and set $\lambda'=\lambda_1 y$, then
\begin{eqnarray*} \|{|{\xi}|^{-b}
\hat{g}(|{\xi}|+\xi_1)}\|_{L_{\xi}^2}^2 &=& \int_{\R_+^n}
\left(\frac{|{\lambda}|^2}{2\lambda_1}\right)^{-2b}|{\hat{g}(\lambda_1)}|^2
\frac{|{\lambda}|^2 }{2\lambda_1^2}\dd \lambda\\
&=& C \int_0^{\infty} \lambda_1^{2b-2}|{\hat{g}(\lambda_1)}|^2
\int_{\R^{n-1}} \left(\lambda_1^2+|{\lambda'}|^2\right)^{1-2b}\dd
\lambda'
\dd \lambda_1 \\
&= & C \int_0^{\infty} \lambda_1^{n-1-2b}|{\hat{g}(\lambda_1)}|^2
\dd
\lambda_1 \int_{\R^{n-1}} \left(1+|{y}|^2\right)^{\frac{1-n}{2}}\dd y\\
&=&\infty \end{eqnarray*} provided that $g\neq 0$ in $\Sc$, which
contradict to \eqref{I}.
\end{prf}

As a complement to the failure of some $r=\infty$ estimate in
\eqref{Sobolev}, we give here a simple but somewhat interesting
result.
\begin{prop}\label{anti}Let $2\leq q<\infty$ and $b=\frac{n}{ 2}-\frac{1}{ q}$, we
have for any $n$,
$$\|{\exp(i t \D)f(x)}\|_{L_x^\infty L_t^q}
\lesssim \|{f}\|_{\dot H^b}$$
\end{prop}
\begin{prf}Let $\alpha=1/2-1/q$, $M$ denote the space of finite
measure and $\mathrm{S}$ denote the usual spherical measure. Note
that if set $x=r\tilde{\omega}, \xi=\lambda
\omega$,\begin{eqnarray*} \|{\F_{\xi}^{-1}(f)(x)}\|_{L^\infty_x}&
\simeq & \|{\int \e^{i \lambda \omega\cdot x}f(\lambda
\omega)\lambda^{n-1}\dd
\lambda\dd\mathrm{S}(\omega)}\|_{L^\infty_x}\\
&\lesssim& \int\|{   \F_r^{-1}(f(\lambda \omega)
\lambda^{n-1})(\omega\cdot x)
}\|_{L^\infty_x} \dd\mathrm{S}(\omega) \\
&\lesssim&  \int\|{
f(\lambda\omega) \lambda^{n-1}}\|_{M_\lambda}\dd\mathrm{S}(\omega)\\
&\lesssim&\|{ f(\lambda\omega)
\lambda^{n-1}}\|_{L_{\mathrm{S}(\omega)}^2M_\lambda}
 \ ,\end{eqnarray*}
then \begin{eqnarray*} \|{\exp(i t \D)f(x)}\|_{L_x^\infty
L_t^q}&\lesssim&\|{\exp(i t \D)f(x)}\|_{L_x^\infty \dot H_t^\alpha}\\
&\simeq&
\|{\F_\xi^{-1}\left(\delta(\tau-|{\xi}|)|{\tau}|^\alpha\hat{f}(\xi)\right)}\|_{L_x^\infty
L_\tau^2}\\
&\lesssim&
\|{\F_\xi^{-1}\left(\delta(\tau-|{\xi}|)|{\tau}|^\alpha\hat{f}(\xi)\right)}\|_{L_\tau^2
L_x^\infty}\\
&\lesssim&
\|{\delta(\tau-\lambda)|{\tau}|^\alpha\hat{f}(\lambda\omega)\lambda^{n-1}}\|_{L_\tau^2L_{\mathrm{S}(\omega)}^2
M_\lambda} \\
&=&\|{|{\tau}|^{n-1+\alpha}\hat{f}(\tau\omega)}\|_{L_\tau^2L_{\mathrm{S}(\omega)}^2}\\
&=&
\|{|{\xi}|^{\frac{n-1}{2}+\alpha}\hat{f}(\xi)}\|_{L_\xi^2}\simeq
\|{f}\|_{\Hs^{b}}
\end{eqnarray*}
\end{prf}

\section{Radial Improvement}
Now we turn to the radial or angular improvements of
\eqref{Sobolev}. Recently, Sterbenz \cite{Stbz04p} gets some
improvement for $n\geq 4$ in condition that the data with some
addition angular regularity. As remarked there, the $n=3$
counterpart follows directly as a combination of result in
\cite{MaNaNaOz03p} with the Proposition 3.4 in \cite{Stbz04p}. We
summarize the complete results here.
\begin{thm}\label{Ster}
Let $n\geq 3$ be the number of spatial dimensions,
$\sigma_\Omega=n-1$, $\sigma=\frac{n-1}{2}$. Then for every
$\epsilon>0$, there is a $C_\epsilon$ depending only on $\epsilon$
such that the following set of estimates hold for any $f\in\Sc$:

\begin{equation} \label{ang} \|{\e^{i t\D}f(x)}\|_{L_t^qL_x^r} \ \lesssim \
C_\epsilon\,\left(
        \|{  \langle\Omega\rangle^s
        f(x)}\|_{\dot H^b} \right)
    \ ,\end{equation}
where we have that $r<\infty$, \ $s= (1 + \epsilon)(\frac{n-1}{r}
+ \frac{2}{q}- \frac{n-1}{2})$, \ $\frac{1}{q} + \frac{n}{r} =
\frac{n}{2} -b$\ , \ $\frac{1}{q} + \frac{\sigma}{r} \geqslant
\frac{\sigma}{2}$\ , and \ \ $\frac{1}{q} +
\frac{\sigma_\Omega}{r} < \frac{\sigma_\Omega}{2}$\ . All of the
implicit constants in the above inequality depend on $n$, $q$, and
$r$. Here $\Omega_{i,j} \ := \ x_i\partial_j - x_j\partial_i $,
$\Delta_{sph} \ := \ \sum_{i<j} \Omega_{ij}^2$, $|\Omega|^s =
(-\Delta_{sph})^\frac{s}{2}$, and
\begin{equation}
    \|{\langle\Omega\rangle^s f}\|_{\dot H^b}
    \ = \ \|{f}\|_{\dot H^b} + \|{|\Omega|^s
    f}\|_{\dot H^b}\ . \notag
\end{equation}
\end{thm}
Note that this result is only proved for the case $n\geq 3$, and
it seems that the argument in \cite{Stbz04p} can only work for
$n\geq 3$. It would be interesting to extend this result to
$n=2$(compare with Theorem \ref{radi}).

Now we give the proof of
Theorem \ref{radi}.\\
{\bf Proof of Theorem \ref{radi}} In view of Theorem \ref{Ster}
and the usual interpolation, the theorem is reduced to the proof
of $r=\infty$ case. Let $|x|=r$, $|\xi|=\lambda$,
$x\cdot\xi=r\lambda\cos \theta=r \lambda y$ and
$\hat{g}(\lambda):=\lambda^{n-1}\hat{f}(\lambda)H(\lambda)\in L^1$
with $H$ the usual Heaviside function. Then \begin{eqnarray*}\e^{i
t\D} f&=&\int \e^{i (x\cdot \xi+t|\xi|)}\hat{f}(\xi)\dd\xi\\&
\simeq&\int_0^\pi \int^\infty_0 \e^{i \lambda(t+ r
\cos \theta)}\hat{f}(\lambda)\lambda^{n-1}(\sin \theta)^{n-2}\dd\lambda\dd\theta\\
&\simeq & \int_0^\pi g(t+ r \cos\theta)(\sin\theta)^{n-2}\dd\theta\\
&=&\int_{-1}^{1}g(ry+t)(1-y^2)^{\frac{n-3}{2}}\dd y:=I(r,t)\
,\end{eqnarray*} Now if $n\geq 3$, $(1-y^2)^{\frac{n-3}{2}}\leq
1$, $|I(r,t)|\leq\int_{-1}^{1}|g(ry+t)|\dd
y=\frac{1}{r}\int_{t-r}^{t+r}|g(z)|\dd z$, then for any $2\leq
q<\infty$,
$$\|\e^{i t\D}
f\|_{L_t^q L_x^\infty} \lesssim \|I\|_{L_t^q L_r^\infty}
\lesssim\|\mathcal{M}g(t)\|_{L^q} \lesssim
 \|g\|_{L^q} \lesssim \|g\|_{\Hs^{1/2-1/q}}\simeq
\|f\|_{\Hs^{n/2-1/q}}\ ,$$ here $\mathcal{M}$ is the usual maximal
operator. If $n=2$, note that $(1-y^2)^{-\frac{1}{2}}\in
L^{p'}_{J}$ with $J=[-1,1]$ and $p\in(2,\infty)$,
$|I(r,t)|\lesssim \|g(t+ry)\|_{L^p_{y\in J}}$, then for any
$q\in(2,\infty)$, let $p=1+q/2$, \begin{eqnarray*} \|\e^{i t\D}
f\|_{L_t^q L_x^\infty}^q&\lesssim&\|g(t+ry)\|^q_{L_t^q
L_r^\infty L^p_{y\in J}}\\
&=&\||g(t+ry)|^p\|^{q/p}_{L_t^{q/p}
L_r^\infty L^1_{y\in J}}\\
&\lesssim&\|\mathcal{M}(|g|^p)(t)\|^{q/p}_{L_t^{q/p}}\\
&\lesssim&\|g\|_{L^q}^q\lesssim\|g\|^q_{\Hs^{1/2-1/q}}\simeq\|f\|_{\Hs^{1-1/q}}^q\end{eqnarray*}
\hfill {\vrule height6pt width6pt depth0pt}\medskip

\section{Inhomogeneous Space $H^s$}
In this section, based on Theorem \ref{gen}, we use scaling
argument and Sobolev embedding to derive estimate in
$H^s=(1-\triangle)^{-s/2}L^2$ instead of $\Hs^b$, i.e., we give
the proof of Theorem \ref{inho} and Theorem \ref{d-inh}. Note that
$u=\cos(t\D)f+\D^{-1}\sin(t\D)g=v+w$ solve the homogeneous wave
equation $\square u=0$ with data $(f,g)$, we'll deal with the two
parts separately
below. In this section, we use $b+$ to denote $b+\epsilon$ with
$\epsilon>0$ arbitrary small.

Note that for any $(q,r,n)$ admissible, we have $b\geq 0$. Then
for all admissible $(q,r,n)$ except that
$(q,r)=(\max(2,\frac{4}{n-1}),\infty)$ and
$(q,r)=(\infty,\infty)$, we have the $H^b$ estimate for
$v=\cos(t\D)f$ and $\tilde{v}=\exp(i t \D)f$. For
$(q,r)=(\infty,\infty)$, the failure of Sobolev embedding
$H^{n/2}\subset L^\infty$ gives the failure of the $H^{n/2}$
estimate. On the other hand, for $(q,\infty,n)$ with
$q>\frac{4}{n-1}$ and $q\geq 2$, there exists $p$ such that
\eqref{Sobolev} valid for $(q,p,n)$,
$$ \|v\|_{L_t^qL_x^\infty}\lesssim \|(1+\D)^{\frac{n}{p}+}v\|_{L_t^q L_x^{p}}\lesssim \|f\|_{H^{b+}}\ .$$
The argument in Proposition \ref{gen-f} yields that one need
$s>\frac{n+1}{4}$ for the $H^s$ estimate of $\tilde{v}$ with
$r=\infty$.
Note that $$H^s\ \mathrm{estimate\ valid\ for\ some}\ s\Rightarrow
\mathrm{single\ frequency\ estimate}$$ then no $H^s$ estimate
valid for $(2,\infty,3)$. By duality argument and decay estimate,
we get the $(4,\infty,2)$ estimate for $f$ with $\hat{f}$ support
in unit ball, then we get $(4,\infty,2)$ estimate with $s=b+$. In
conclusion, we have proved and the following(both see Figure
\ref{drv-p})

\begin{prop}\label{1st} Let $n\geq 2$. We have \begin{equation} \label{shi}\|\cos
(t\D)f\|_{L_t^q L_x^r}\lesssim \|f\|_{H^s}\end{equation} with
$s=b+$ iff $(q,r,n)$ admissible and $(q,r,n)\neq (2,\infty,3)$.
Moreover, for admissible $(q,r,n)$ except that
$(q,r)=(\max(2,\frac{4}{n-1}),\infty)$ and
$(q,r)=(\infty,\infty)$, we have \eqref{shi} valid with $s=b$. On
the other hand, if \eqref{shi} valid with $s=b$, then $(q,r,n)$ is
admissible and $(q,r,n)\neq (2,\infty,3),(\infty,\infty,n)$.
\end{prop}
\begin{prop}\label{deriv} The same result as in Proposition \ref{1st} valid for the estimate
\begin{equation}\label{derived}\|\exp (i t\D)f\|_{L_t^q L_x^r}\lesssim
\|f\|_{H^s}\end{equation} Moreover, \eqref{derived} fail to be
valid if $(q,r)=(\frac{4}{n-1},\infty)$ and $s=b$.
\end{prop}
And hence we have the Theorem \ref{d-inh}.

For $w=\D^{-1}\sin(t\D)g$, the situation is different. A scaling
argument yields
\begin{prop}\label{2nd}If
$b<1$, then for any $s\in\mathbb{R}$ the estimate
\begin{equation} \label{shi2}\|w\|_{L_t^q L_x^r}\lesssim \|g\|_{H^{s-1}}\end{equation}
fails. If $b=1$, then \eqref{shi2} valid for $s=b+$ if and only if
\eqref{shi2} with $s=b$ valid.
\end{prop}
\begin{prf}
For the $b<1$ part, we need only to show the failure for large
$s$. For such $s$, let
$w_{\lambda}(t,x)=\lambda^{\frac{1}{q}+\frac{n}{r}}w(\lambda
t,\lambda x)$, and hence
$g_{\lambda}(x)=\lambda^{\frac{1}{q}+\frac{n}{r}+1}g(\lambda x)$.
Then if \eqref{shi2} valids, we have for $s\geq 1$,
$$\begin{array}{l}\|w\|_{L_t^qL_x^r}=\|w_\lambda\|_{L_t^qL_x^r}\lesssim\|g_\lambda\|_{H^{s-1}}\lesssim
\|g_\lambda\|_{\Hs^{s-1}}+\|g_\lambda\|_{L^2}\\
=\lambda^{s-b}\|g\|_{\Hs^{s-1}}+\lambda^{1-b}\|g\|_{L^2}\rightarrow
0
\end{array}$$ by letting $\lambda\rightarrow 0$, so is contradicted.

For $b=1$, if \eqref{shi2} valid for $s=b+\epsilon$, then
$$\begin{array}{l}\|w\|_{L_t^q
L_x^r}=\|w_\lambda\|_{L_t^qL_x^r}\lesssim\|g_\lambda\|_{H^{\epsilon}}\lesssim
\|g_\lambda\|_{\Hs^{\epsilon}}+\|g_\lambda\|_{L^2}\\
=\lambda^{\epsilon}\|g\|_{\Hs^{\epsilon}}+\|g\|_{L^2}\rightarrow
\|g\|_{L^2}\end{array}$$ with  $\lambda\rightarrow 0$.
\end{prf}

\begin{figure}\centering
  \includegraphics[width=0.90\textwidth]{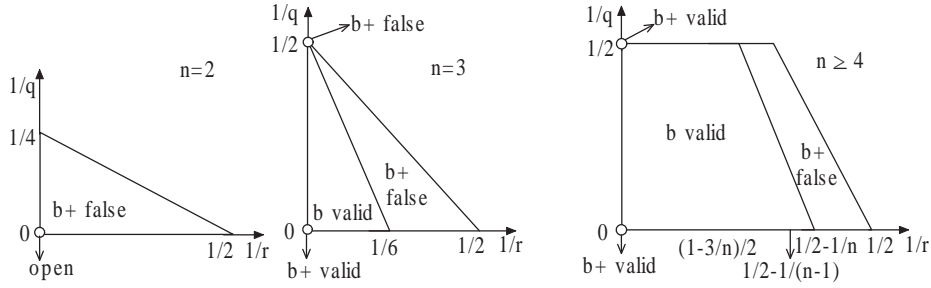}\\
  \caption{Strichartz Estimate \eqref{shi2}}\label{ginh-p}
\end{figure}

If $b>1$, then we have $H^{b-1}$ estimate once we have $\Hs^{b-1}$
estimate.
On the other hand, for $(q,\infty,n)$ with $q>\frac{4}{n-1}$ and
$q\geq 2$ such that $b>1$, there exists $p$ such that
\eqref{Sobolev} valid for $(q,p,n)$ with
$\lambda_1=\frac{n}{2}-\frac{n}{p}-\frac{1}{q}-1\geq 0$,
$$ \|w\|_{L_t^q L_x^\infty}\lesssim \|(1+\D)^{\frac{n}{p}+}w\|_{L_t^q L_x^{p}}\lesssim
\|(1+\D)^{\frac{n}{p}+}\D^{\lambda_1}g\|_{L^2}\lesssim\|g\|_{H^{-1+b+}}\
.$$ Combining Proposition \ref{2nd} with the previous observation,
we get(for its figure, see Figure \ref{ginh-p})
\begin{prop}\label{2ndg} The estimate \eqref{shi2} valid with $s=b+$
if $(q,r,n)$ admissible, $b\geq 1$, $n\geq 3$ and $(q,r,n)\neq
(2,\infty,3)$. On the other hand, $b< 1$ or $(q,r,n)=(2,\infty,3)$
imply that \eqref{shi2} fail with $s=b+$. If $(q,r,n)$ admissible,
$b\geq 1$, $n\geq 3$ and $(q,r)\neq (2,\infty),(\infty,\infty)$,
we have \eqref{shi2} valid with $s=b$.
\end{prop}

\end{document}